\documentclass[reqno,12pt,a4paper]{amsart}

\usepackage{mathrsfs}
\usepackage{amssymb}
\usepackage{amsfonts}
\usepackage{latexsym}
\usepackage{amsthm}

\usepackage{graphicx}
\def\lb{\label}

\newcommand{\er}[1]{\textrm{(\ref{#1})}}

\begin{document}


\renewcommand{\theequation}{\arabic{section}.\arabic{equation}}
\theoremstyle{plain}
\newtheorem{theorem}{\bf Theorem}[section]
\newtheorem{lemma}[theorem]{\bf Lemma}
\newtheorem{corollary}[theorem]{\bf Corollary}
\newtheorem{proposition}[theorem]{\bf Proposition}
\newtheorem{definition}[theorem]{\bf Definition}
\newtheorem{remark}[theorem]{\it Remark}

\def\a{\alpha}  \def\cA{{\mathcal A}}     \def\bA{{\bf A}}  \def\mA{{\mathscr A}}
\def\b{\beta}   \def\cB{{\mathcal B}}     \def\bB{{\bf B}}  \def\mB{{\mathscr B}}
\def\g{\gamma}  \def\cC{{\mathcal C}}     \def\bC{{\bf C}}  \def\mC{{\mathscr C}}
\def\G{\Gamma}  \def\cD{{\mathcal D}}     \def\bD{{\bf D}}  \def\mD{{\mathscr D}}
\def\d{\delta}  \def\cE{{\mathcal E}}     \def\bE{{\bf E}}  \def\mE{{\mathscr E}}
\def\D{\Delta}  \def\cF{{\mathcal F}}     \def\bF{{\bf F}}  \def\mF{{\mathscr F}}
\def\c{\chi}    \def\cG{{\mathcal G}}     \def\bG{{\bf G}}  \def\mG{{\mathscr G}}
\def\z{\zeta}   \def\cH{{\mathcal H}}     \def\bH{{\bf H}}  \def\mH{{\mathscr H}}
\def\e{\eta}    \def\cI{{\mathcal I}}     \def\bI{{\bf I}}  \def\mI{{\mathscr I}}
\def\p{\psi}    \def\cJ{{\mathcal J}}     \def\bJ{{\bf J}}  \def\mJ{{\mathscr J}}
\def\vT{\Theta} \def\cK{{\mathcal K}}     \def\bK{{\bf K}}  \def\mK{{\mathscr K}}
\def\k{\kappa}  \def\cL{{\mathcal L}}     \def\bL{{\bf L}}  \def\mL{{\mathscr L}}
\def\l{\lambda} \def\cM{{\mathcal M}}     \def\bM{{\bf M}}  \def\mM{{\mathscr M}}
\def\L{\Lambda} \def\cN{{\mathcal N}}     \def\bN{{\bf N}}  \def\mN{{\mathscr N}}
\def\m{\mu}     \def\cO{{\mathcal O}}     \def\bO{{\bf O}}  \def\mO{{\mathscr O}}
\def\n{\nu}     \def\cP{{\mathcal P}}     \def\bP{{\bf P}}  \def\mP{{\mathscr P}}
\def\r{\rho}    \def\cQ{{\mathcal Q}}     \def\bQ{{\bf Q}}  \def\mQ{{\mathscr Q}}
\def\s{\sigma}  \def\cR{{\mathcal R}}     \def\bR{{\bf R}}  \def\mR{{\mathscr R}}
\def\S{\Sigma}  \def\cS{{\mathcal S}}     \def\bS{{\bf S}}  \def\mS{{\mathscr S}}
\def\t{\tau}    \def\cT{{\mathcal T}}     \def\bT{{\bf T}}  \def\mT{{\mathscr T}}
\def\f{\phi}    \def\cU{{\mathcal U}}     \def\bU{{\bf U}}  \def\mU{{\mathscr U}}
\def\F{\Phi}    \def\cV{{\mathcal V}}     \def\bV{{\bf V}}  \def\mV{{\mathscr V}}
\def\P{\Psi}    \def\cW{{\mathcal W}}     \def\bW{{\bf W}}  \def\mW{{\mathscr W}}
\def\o{\omega}  \def\cX{{\mathcal X}}     \def\bX{{\bf X}}  \def\mX{{\mathscr X}}
\def\x{\xi}     \def\cY{{\mathcal Y}}     \def\bY{{\bf Y}}  \def\mY{{\mathscr Y}}
\def\X{\Xi}     \def\cZ{{\mathcal Z}}     \def\bZ{{\bf Z}}  \def\mZ{{\mathscr Z}}
\def\O{\Omega}

\newcommand{\gA}{\mathfrak{A}}
\newcommand{\gB}{\mathfrak{B}}
\newcommand{\gC}{\mathfrak{C}}
\newcommand{\gD}{\mathfrak{D}}
\newcommand{\gE}{\mathfrak{E}}
\newcommand{\gF}{\mathfrak{F}}
\newcommand{\gG}{\mathfrak{G}}
\newcommand{\gH}{\mathfrak{H}}
\newcommand{\gI}{\mathfrak{I}}
\newcommand{\gJ}{\mathfrak{J}}
\newcommand{\gK}{\mathfrak{K}}
\newcommand{\gL}{\mathfrak{L}}
\newcommand{\gM}{\mathfrak{M}}
\newcommand{\gN}{\mathfrak{N}}
\newcommand{\gO}{\mathfrak{O}}
\newcommand{\gP}{\mathfrak{P}}
\newcommand{\gQ}{\mathfrak{Q}}
\newcommand{\gR}{\mathfrak{R}}
\newcommand{\gS}{\mathfrak{S}}
\newcommand{\gT}{\mathfrak{T}}
\newcommand{\gU}{\mathfrak{U}}
\newcommand{\gV}{\mathfrak{V}}
\newcommand{\gW}{\mathfrak{W}}
\newcommand{\gX}{\mathfrak{X}}
\newcommand{\gY}{\mathfrak{Y}}
\newcommand{\gZ}{\mathfrak{Z}}

\def\ve{\varepsilon}   \def\vt{\vartheta}    \def\vp{\varphi}    \def\vk{\varkappa}

\def\Z{{\mathbb Z}}    \def\R{{\mathbb R}}   \def\C{{\mathbb C}}    \def\K{{\mathbb K}}
\def\T{{\mathbb T}}    \def\N{{\mathbb N}}   \def\dD{{\mathbb D}}


\def\la{\leftarrow}              \def\ra{\rightarrow}            \def\Ra{\Rightarrow}
\def\ua{\uparrow}                \def\da{\downarrow}
\def\lra{\leftrightarrow}        \def\Lra{\Leftrightarrow}


\def\lt{\biggl}                  \def\rt{\biggr}
\def\ol{\overline}               \def\wt{\widetilde}
\def\no{\noindent}


\let\ge\geqslant                 \let\le\leqslant
\def\lan{\langle}                \def\ran{\rangle}
\def\/{\over}                    \def\iy{\infty}
\def\sm{\setminus}               \def\es{\emptyset}
\def\ss{\subset}                 \def\ts{\times}
\def\pa{\partial}                \def\os{\oplus}
\def\om{\ominus}                 \def\ev{\equiv}
\def\iint{\int\!\!\!\int}        \def\iintt{\mathop{\int\!\!\int\!\!\dots\!\!\int}\limits}
\def\el2{\ell^{\,2}}             \def\1{1\!\!1}
\def\sh{\sharp}
\def\wh{\widehat}
\def\bs{\backslash}
\def\na{\nabla}

\def\sh{\mathop{\mathrm{sh}}\nolimits}
\def\all{\mathop{\mathrm{all}}\nolimits}
\def\Area{\mathop{\mathrm{Area}}\nolimits}
\def\arg{\mathop{\mathrm{arg}}\nolimits}
\def\const{\mathop{\mathrm{const}}\nolimits}
\def\det{\mathop{\mathrm{det}}\nolimits}
\def\diag{\mathop{\mathrm{diag}}\nolimits}
\def\diam{\mathop{\mathrm{diam}}\nolimits}
\def\dim{\mathop{\mathrm{dim}}\nolimits}
\def\dist{\mathop{\mathrm{dist}}\nolimits}
\def\Im{\mathop{\mathrm{Im}}\nolimits}
\def\Iso{\mathop{\mathrm{Iso}}\nolimits}
\def\Ker{\mathop{\mathrm{Ker}}\nolimits}
\def\Lip{\mathop{\mathrm{Lip}}\nolimits}
\def\rank{\mathop{\mathrm{rank}}\limits}
\def\Ran{\mathop{\mathrm{Ran}}\nolimits}
\def\Re{\mathop{\mathrm{Re}}\nolimits}
\def\Res{\mathop{\mathrm{Res}}\nolimits}
\def\res{\mathop{\mathrm{res}}\limits}
\def\sign{\mathop{\mathrm{sign}}\nolimits}
\def\span{\mathop{\mathrm{span}}\nolimits}
\def\supp{\mathop{\mathrm{supp}}\nolimits}
\def\Tr{\mathop{\mathrm{Tr}}\nolimits}
\def\BBox{\hspace{1mm}\vrule height6pt width5.5pt depth0pt \hspace{6pt}}
\def\where{\mathop{\mathrm{where}}\nolimits}
\def\as{\mathop{\mathrm{as}}\nolimits}


\newcommand\nh[2]{\widehat{#1}\vphantom{#1}^{(#2)}}
\def\dia{\diamond}

\def\Oplus{\bigoplus\nolimits}



\def\qqq{\qquad}
\def\qq{\quad}
\let\ge\geqslant
\let\le\leqslant
\let\geq\geqslant
\let\leq\leqslant
\newcommand{\ca}{\begin{cases}}
\newcommand{\ac}{\end{cases}}
\newcommand{\ma}{\begin{pmatrix}}
\newcommand{\am}{\end{pmatrix}}
\renewcommand{\[}{\begin{equation}}
\renewcommand{\]}{\end{equation}}
\def\eq{\begin{equation}}
\def\qe{\end{equation}}
\def\[{\begin{equation}}
\def\bu{\bullet}

\title[{Hamiltonian and small  action variables for periodic dNLS  }]
{Hamiltonian and small  action variables  for periodic dNLS}

\date{\today}
\address{School of Math., Cardiff University.
Senghennydd Road, CF24 4AG Cardiff, Wales, UK.
email \ KorotyaevE@cf.ac.uk,
{\rm Partially supported by EPSRC grant EP/D054621.}}
\author[Evgeny L. Korotyaev]{Evgeny L. Korotyaev}



\begin{abstract}
We consider the defocussing  NLS equation with small periodic initial condition.
A new approach to study the Hamiltonian as a function of action variables is demonstrated.
The problems for the NLS equation is reformulated  as the problem of
conformal mapping theory corresponding to quasimomentum of the Zakharov-Shabat operator.
The main tool is the L\"owner type equation for the quasimomentum.
In particular, we determine the asymptotics of the Hamiltonian for small action variables.
Moreover, we determine the gradient of Hamiltonian with respect to action variables.
This gives so called frequencies and determines how the angles variables depend on the time.

\end{abstract}

\maketitle

{\it Dedicated to the memory of my teacher Mikhail Sh. Birman, 1928-2009}

\section{\bf Introduction}

Consider the defocussing cubic non-linear Schr\"odinger equation (dNLS)
$$
J{\pa\p \/\pa t}=-\p _{xx}+2|\p |^2\p ,\ \ \ J=\ma 0&1\\-1&0\am,\qqq
\p =\ma \p_1\\ \p_2\am ,\qq
$$
on the circle $\T=\R/\Z$, i.e. $\p(x + 1,t) = \p(x, t)$ for $x, t\in \R$, with
the initial conditions:
$$
\p(\cdot,0)=q=\ma q_1\\ q_2\am\in L^2(\T)\oplus L^2(\T),\qq
$$
where $q_1, q_2$ are real functions.
 The dNLS equation has the Hamiltonian $H$ given by
\[
\lb{Ha}
H(q)={1\/2}\int_0^1(|q'(x)|^2+|q(x)|^4)dx.
\]
For functionals $E=E(q)$ and $G=G(q)$ the Poisson bracket $\{E,G\}_P$ has  the form
$$
\{E,G\}_P=\int _0^1\rt({\pa E \/  \pa q_1(x)} {\pa G \/\pa q_2(x)}
-{\pa G\/\pa q_1(x)} {\pa E\/\pa q_2(x)}\rt)dx.
$$
The Hamiltonian system with Hamiltonian $H$ is given by
$$
{\pa\p \/\pa t}=\{H(\p),\p\}_P=-J \ma {\pa H\/\pa \p_1}\\{\pa H\/\pa\p_2}\am  =-J(-\p_{xx}+2|\p |^2\p).
$$

 The dNLS equation has the Hamiltonian $H$ and other two integrals $H_0$ and $H_1$ given by
$$
H_0(q)=\|q\|^2=\int_0^1(q_1(x)^2+q_2(x)^2)dx,
$$$$
H_1(q)=\int_0^1(q_2'(x)q_1(x)-q_1'(x)q_2(x))dx.
$$
 There are many papers and books devoted to the dNLS equation, see \cite{FT},
 \cite{GH}.
  The action-angle variables for the dNLS equation  on the circle were studied by
McKean-Vaninsky \cite{MV1}, \cite{MV2}, and
  Grebert-Kappeler-P\"oschel \cite{GKP}, see also \cite{V}.
  Note that the action-angle variables are studied in \cite{VN}, \cite{Ku}, \cite{KP}.
The dNLS equation admits globally defined real analytic
action-angle variables $A_n, \f_n, n\in \Z$ (see \cite{GKP}, \cite{MV1}, \cite{MV2}, \cite{VN}).
Recall the important identities from \cite{K2}:
\[
\lb{A0}
H_0=\|q\|^2=\sum_{n\in\Z} A_n,
\]
\[
\lb{A1}
H_1=\sum_{n\in\Z} (2\pi n)A_n,
\]
\[
\lb{A22}
H=\sum_{n\in\Z} (2\pi n)^2A_n+2H_0^2-U,\qq
\]
where $U$ is some nonlinear functional given by \er{AV}.
In fact we rewrite $H_0, H_1$
and the main part of $H$ in terms of simple functions of the actions $A=(A_n)_{n\in \Z}$.
Moreover, the following estimates from \cite{K2} hold true:
\[
\lb{eV}
0\le U(A)\le {4\/3}\|A\|_1^2,\qqq \qqq if \qq |q|\in L^2(\T),
\]
where $\|A\|_1=\sum_{n\in\Z} A_n$.
Introduce the real  spaces
$$
 \ell^p=\{f=(f_n)_{n\in \Z},\ \ \| f \|_p<\iy \},
\qqq \ \ \| f \|_p^p=\sum _{n} f_n^p <\iy, \ p\ge 1.
$$
Note that \er{A0}, \er{A22}, \er{eV}, yield:

\no (i) $|q|\in L^2(\T)$ iff $(A_n)_{n\in \Z}\in\ell^1$,

\no (ii) $|q'|\in L^2(\T)$ iff $(nA_n)_{n\in \Z}\in\ell^1$.

In Theorem \ref{T1} we will show that  $U(A)$ is a well defined function of $A=(A_n)_{n\in \Z}\in \ell^2$.
We mention the Kuksin conjecture \cite{Ku1}:
{\it there is an estimate of the Hamiltonian $H$ in terms of $\cP(\|A\|_2)$ for some polynomial $\cP(z)$}.
As far as we know, no an estimate of $H$ in terms of $\cF(\|A\|_2)$ for some function $\cF(z), z\ge 0$
has been published.

We formulate our result about the estimates.

\begin{theorem}
\lb{T1}
The following estimates hold true:
\[
\lb{eV1}
{\pi\/6}\|A\|_2^2\le U(A)\le {2\pi\/3}\sqrt{C_1}\|A\|_2^2,\qqq if \qqq  (|\hat q_n|)_{n\in\Z}\in \ell^4,
\qqq
\]
where  $\hat q_n=\int_0^1q(x)e^{-i2\pi nx}dx$
and $C_1=\max \{2, \cosh {\pi\/2}\|A\|_\iy\}$.
\end{theorem}

\no {\bf Remark.} 1) $A\in\ell^2$ iff $(|\hat q_n|)_{n\in\Z}\in \ell^4$, see Sect. 2.

\no 2) If $q\in L^{4\/3}(\T)$, then $(|\hat q_n|)_{n\in\Z}\in \ell^4$.

\no 2) \er{eV1} gives $U(A)=0$ for some $A\in\ell^2$ iff $A=0$ (or $q=0$).

The Hamiltonian $H$ depends only on the actions $A$.
Introduce the frequencies $ \O_n$ by
\[
\O_n=\pa_n H, \qq n\in \Z,\qqq \pa_n={\pa \/\pa A_n}, \qqq \pa =(\pa_n)_{n\in \Z}.
\]
The parameters $\O_n$ are very important, since the angle variables $\f_n(t)$ as functions of time $t\ge 0$  have the form
$$
\f_n(t)=\f_n(0)+\O_n t, \qq t\ge 0, \qq n\in \Z.
$$
Due to \er{A2} we deduce that the gradient is given by
\[
\O_n=\pa_n H=(2\pi n)^2+4H_0-\pa_n U. \qqq 
\]
Thus in order to study $\O_n$ we need to study $\pa_n U$ only, which is defined on $q\in L^2(\T)$.  Our goal is to give a new method to study
Hamiltonian as a function of action variables.
In this paper we reformulate the problems for the dNLS equation as the problems of the conformal mapping theory.
The main technical tool
is Theorem \ref{Tdek} from \cite{KK3} about the L\"owner type equation. We formulate our main result.

\begin{theorem}
\lb{T2}
Let $\mA=\{A\in \ell^1: \sum_{n\in \Z}A_n\le  {1\/8^2}\}$.
 Then the function $U: \mA\to [0,\iy)$ has the derivative ${\pa_n} U(A)$
for each $n\in \Z$, which is continuous on $\mA$ and satisfies
\[
\lb{TV}
|U(A)-\|A\|_2^2|\le 4\pi \sqrt 3 \|A\|_2^3,
\]
\[
\lb{To}
\|\pa U(A)-2A\|_2\le 11\pi^2\|A\|_\iy \|A\|_2.
\]
\end{theorem}

\no {\bf Remark.} 1) Note that the Hamiltonian $H(q)$ is defined on the functions $|q'|\in L^2(\T)$, but we study the frequencies for the potentials $|q|\in L^2(\T)$.

\no 2) In order to prove \er{TV} with the reminder of the third order $\|A\|_2^3$
we use the properties of the quasimomentum from \cite{KK2}, \cite{K3},\cite{K4} and the identity \er{A2} only.
Moreover, similar arguments give the Taylor series of the Hamiltonian $H$ in terms of actions $A_n$ for any order.
We consider the simplest case to write the short paper and to formulate our approach. The main problem is to show \er{To}.

\no 3) We estimate $U$ in terms of the norm $\|A\|_2$,
which yields estimates in terms of $\|A\|_1$, since
$\|A\|_2\le \|A\|_1$. In order to study the perturbation of the dNLS the estimates the Hamiltonian $H(q)$ in terms of $\|A\|_2$ are important.
This is the motivation of the estimates \er{TV}, \er{To}.

\no 4)  In the case  $|q'|\in L^2(\T)$ the Hamiltonian $H$
is a real analytic function of $(nA_n)_{n\in \Z}$
and the Marchenko-Ostrovski parameters, similar to the case
 $H_0$, see Theorem \ref{Tdek}.

\no 4) In \cite{K5} we consider the frequencies $\O_n$ for general case and determine their asymptotics  as $n\to \pm \iy$.

We now describe the plan of the paper.  In Section 2 we recall the needed results
from \cite{KK2}, \cite{KK3}, \cite{K4}, \cite{K3} and \cite{MO} about the quasimomentum $k(z,h)$
as the function of two variables $z,h$:  $z$ is a spectral variable and $h\in \ell^2$ is the Marchenko-Ostrovski parameter. Note that for each $h\in \ell^2 $ the function  $k(\cdot,h):\cK(h)\to \cZ$ is the conformal mapping  (see definitions of $\cK(h), \cZ$ and $k(\cdot,h)$ in \er{dequ}).
 In Theorem \ref{Tdek} from \cite{KK3}  we recall the needed properties of the function $z(k,h)$ (for fix $h\in \ell^2$ the function $z(k,h)$ is an inverse for $k(z,h)$ ) as the function of quasimomentum $k$ and $h\in \ell^2$ and formulate the L\"owner type equation \er{Tdek-2} for the conformal mapping $z(\cdot,h):\cK(h)\to \cZ$.
In Lemma \ref{T7.6} from \cite{KK3} we describe the properties of the actions $A_n(h)$ as functions of the Marchenko-Ostrovski parameters $h\in \ell^2$
and present the exact formulas for ${\pa \/\pa h_m} A_n(h), n,m\in \Z$.
In the proof of Theorem \ref{T1} we use the identity $\o_m=\sum_{n\in \Z} X_{m,n}\wt\o_n$, where $\wt\o_m ={\pa U\/\pa h_m^2}$ and $X_{m,n}={\pa h_m^2 \/\pa A_n}, m\in \Z$. In Lemma \ref{TdA}
we show that operator $X:\ell^2\to \ell^2$ with the matrix $X_{m,n}$
satisfies $X=I+o(1)$ as $A\to 0$.
In Lemma \ref{TgH} we determine $\wt\o_n=2A_n+o(\|A\|)$
as $A\to 0$. Then roughly speaking we determine asymptotics $\o_n={\pa U\/\pa A_n}$.

\section {\bf Preliminaries }
\setcounter{equation}{0}

The dNLS equation is integrable and admits a Lax-pair formalism, see \cite{ZS}, \cite{FT}.
Consider the corresponding Zakharov-Shabat operator $T_{zs}$ acting in $L^2(\R )\oplus L^2(\R )$ and given by
\[
 T_{zs}=J{d\/dx}+\cQ, \qqq \cQ=\ma q_1&q_2\\q_2&-q_1\am, \qq
 \]
where $|q|\in L^1(0,1)$. This operator is essentially self-adjoint
on the domain $\mD=\{f,f',f''\in L^2(\R)\os L^2(\R)\}$, (see \cite{LM}).
  The spectrum of  $T_{zs}$ is purely absolutely continuous and is union of spectral bands $\s_n, n\in \Z$, where
$$
 \s _n=[z_{n-1}^+,z_n^-], \ \dots <z_{2n-1}^-\leq z_{2n-1}^+<z_{2n}^- \leq z_{2n}^+< \dots ,\qq  and
\qq  \ z_n^{\pm}=n\pi +o(1)\    \as  |n|\to\iy .
$$
 The intervals $\s_n$ and $\s_{n+1}$ are separated by gap $g_n=(z^-_n, z^+_n )$ with the length $|g_n|\ge 0$.
 If a gap $g_n$ is degenerate, i.e., $|g_n|=0,$  then the corresponding segments $\s_n, \s_{n+1}$ merge. We use the Zakharov-Shabat  equation for a vector -function $f$:
\[
\lb{zse1}
Jf'+\cQ f=zf, \ \ \ z\in \C,\ \ \ \ \ \ f=\ma f_1\\f_2\am ,
\]
where $f_1, f_2$ are  functions of $x\in\R $. Here and below
$(\  ')=\pa /\pa x $.
The boundary value problem with Eq. \er{zse1} and   with the condition $f(0)=f(1)$
is called periodic and the boundary value problem \er{zse1} with the condition $f(0)=-f(1)$
is called antiperiodic. Here $z_{2n}^{\pm}, n\in \Z$ are the eigenvalues of
the periodic problem  and $z_{2n+1}^{\pm}, n\in \Z$ are the eigenvalues of the
anti-periodic problem. Define the $2\ts 2$-matrix valued fundamental solution  $\P=\P(x,z)$ by
\[
J{d\/dx}\P+\cQ \P=z \P, \ \  \P(0,z)=\ma 1& 0\\ 0&1\am , \ \ z\in \C.
\]
Introduce the Lyapunov function $\D(z)$ by
$$
\D(z)={1 \/ 2}{\rm Tr}  \P(1,z), \ \ z\in \C.
$$
The function $\D$ is entire and $\D(z_n^{\pm})=(-1)^n$ for all $n\in\Z$ (see e.g. \cite{LS}).

We recall results which is crucial for the present paper. For each $q\in L^1(\T)$ there exists a unique
conformal mapping (the quasimomentum) $k:\cZ\to \cK(h) $ with asymptotics
$k(z)=z+o(1)$ as $|z|\to \iy$ (see Fig. 1 and 2) and such that (see  \cite{MO} and   \cite{Mi1}, \cite{Mi2}, \cite{KK2})
\begin{multline}
\lb{dequ}
\cos k(z)= \D(z),\ \  z\in \cZ =\C\sm\cup \ol g_n,\qq and
\ \ \ \
\\
\cK(h)=\C\sm\cup \G_n ,\ \  \  \G_n=(\pi n-i|h_n|,\pi n+i|h_n|), \qqq h_n\in \R,
\\
|h_n|\ge 0 \qqq {\rm is \ defined \ by  \ the \ equation}\qqq  \cosh |h_n| = (-1)^n\D(z_n) \ge 1.
\end{multline}
Here $\G_n$ is the vertical cut and recall that $z_n\in [z_n^-,z_n^+]$ and $\D'(z_n)=0$.
Moreover, we have $(|h_n|)_{n\in \Z}\in \ell^2$ iff $q\in L^2(0,1)$ (and
$(n|h_n|)_{n\in \Z}\in \ell^2$ iff $q'\in L^2(\T)$), see \cite{K1}, \cite{K2}.

Recall the properties of the conformal mapping $k=u(z)+iv(z), z\in \cZ$ from \cite{MO} or \cite{KK2}:

\begin{lemma}
\lb{pk}
Let $h\in \ell^\iy$. Then the quasimomentum $k=u(z)+iv(z), z\in \cZ$
satisfies:

\no 1) $v(z)\ge \Im z>0$ and $v(z)=-v(\ol z)$ for all $z\in \C_+=\{\Im z>0\}$.

\no 2) $v(z)=0$ for all $z\in \s_n, n\in \Z$.

\no 3)  If some $g_n\ne \es, n\in \Z$, then the function
$v(z+i0)>0$ for all $z\in g_n$, and $v(z+i0)$ has a maximum at $z_n\in g_n$
such that $v(z_n+i0)=|h_n|, v'(z_n)=0$, see Fig. 3,  and $\D'(z_n)=0$
 and
\[
\lb{prq}
v(z+i0)=-v(z-i0)>0, \qqq {v'(z+i0)\/z_n-z}>0,\qqq v''(z+i0)<0,   \qqq \all \  z\in g_n\ne \es,
\]
\[
\lb{em-2}
|g_n|\le 2|h_n|.
\]

\no 4) $u'(z)>0$ on  all $(z_{n-1}^+,z_n^-)$ and  $u(z)=\pi n$ for all $z\in g_n\ne \es, n\in \Z$.

\no 5) The function $k(z)$ maps a horizontal cut (a "gap" ) $[z_n^-,z_n^+]$  onto vertical cut $\G_n$  and
a spectral band $\s_n$ onto the segment $[\pi (n-1), \pi n]$ for all $n\in\Z$.
\end{lemma}

We emphasize that the introduction of the quasi-momentum $k(\cdot)$ provides a natural
labeling of all gaps $g_n$ (including the empty ones!) by
demanding that $k(\cdot )$ maps the cut (a "gap" ) $[z_0^-,z_0^+]$ on the vertical
cut $\G_0=(-ih_0,ih_0)$. This determination of a fixed reference
point will be of important later on. Let $z(\cdot )=k^{-1}:\cK(h)\to \cZ$ be the inverse mapping for $k:\cZ\to \cK(h)$.
Below we  will sometimes write $z(k,h), ,..$, instead of $z(k),,..$, when several $h$ are being dealt with.
Note, that for fixed $k\in\cK(h) $ the function $z(k,h),\ h\in \ell^2_{\R}$ is {\bf even with respect
to each variable } $h_n\in \R, n\in \Z$. We use $h_n\in \R$, since it is more convenient for us than $h_n\in [0,\iy)$.

For any
nondegenerate gap $g_n$ the function $k(\cdot,h) $ has an analytic continuation (from above or from below) across the interval $g_n$. This suffices to extend the function $-i(k(\cdot,h)-\pi n) $ by the symmetry.
Similarly the function $z(\cdot,h) $ has an analytic continuation  (from left or from right) across the
vertical cut  $(\pi n-i|h_n|,\pi n+i|h_n|) $ by the symmetry.

In spirit, such result goes back to the
classical Hilbert Theorem (for a finite number of cuts, see e.g.
\cite{J}) in the conformal mapping theory. A similar theorem for the
Hill operator is technically more complicated (there is a infinite
number of cuts) and was proved by Marchenko-Ostrovski \cite{MO}. The
proof of Misura for the Zakharov-Shabat system follows the general
idea from \cite{MO}. For additional properties of the conformal mapping we also refer to
our previous papers \cite{KK2}, \cite{KK4}, \cite{K2}.

A lot of papers are devoted to the inverse problems for the operator $T_{zs}$.
 Misura \cite{Mi1}, \cite{Mi2} extended the results of \cite{MO} to the periodic Zakharov-Shabat operator (in terms of heights $h=(h_n)_{n\in \Z}\in \ell^2$ of  vertical cuts $\G_n, n\in \Z$, see Sect.2). The author \cite{K1} reproved the results of Misura \cite{Mi1}, \cite{Mi2} by the direct method \cite{KK1}.
 The gap length mappings were considered in \cite{BGGK}, \cite{GG},  \cite{K2}.
 A priori estimates (two sided) of $H_0, H$ in terms of actions, gap lengths,
 Marchenko-Ostrovski parameters, etc. were obtained in \cite{KK2}, \cite{KK3}, \cite{K1}-\cite{K4}.
There are papers devoted to the integrals $H, H_0, H_1$,
see  \cite{KK2}, \cite{KK3}, \cite{K1}-\cite{K4}, where various both identities and estimates of integrals in terms of gap lengths of the Zakharov-Shabat operator, actions variables, the Marchenko-Ostrovski parameters (so-called heights $h_n$, see Sect. 2)
were obtained.

Note that  if $q\in L^2(\T)$, then $v\in L^1(\T)\cap L^\iy(\R)$, see \er{HQ}.

The quasimomentum $k(\cdot)$ has
has asymptotics
$$
k(z)=z-{Q_0+o(1)\/z}\qq if \ \ q\in L^2(\T) \qq  and  \qq k(z)=z-{Q_0\/z}-{Q_1\/z^2}-{Q_2+o(1)\/z^3}\qq if \ \
q'\in L^2(\T)
$$
as $z\to i\iy$, see \cite{K2},
where the functionals $Q_j(h), (n^{j\/2}h_n)_{n\in \Z}\in \ell^2, j=0,1,2$ are given by
\[
\lb{dQ}
Q_j(h)={1\/\pi}\int_\R z^jv(z+i0,h)dz\ge 0,\qqq j=0,1,2,\qq k(z)=u(z)+iv(z).
\]
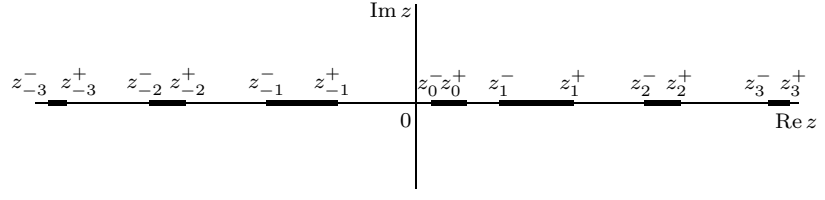
\begin{figure}
\tiny
\unitlength=1mm
\special{em:linewidth 0.4pt}
\linethickness{0.4pt}
\begin{picture}(120.67,34.33)
\put(20.33,21.33){\line(1,0){100.33}}
\put(70.33,10.00){\line(0,1){24.33}}
\put(69.00,19.00){\makebox(0,0)[cc]{$0$}}
\put(120.33,19.00){\makebox(0,0)[cc]{$\Re z$}}
\put(67.00,33.67){\makebox(0,0)[cc]{$\Im z$}}
\put(81.33,21.33){\linethickness{2.0pt}\line(1,0){9.67}}
\put(100.33,21.33){\linethickness{2.0pt}\line(1,0){4.67}}
\put(116.67,21.33){\linethickness{2.0pt}\line(1,0){2.67}}
\put(72.33,21.33){\linethickness{2.0pt}\line(1,0){4.67}}
\put(60.00,21.33){\linethickness{2.0pt}\line(-1,0){9.33}}
\put(40.00,21.33){\linethickness{2.0pt}\line(-1,0){4.67}}
\put(24.33,21.33){\linethickness{2.0pt}\line(-1,0){2.33}}
\put(81.67,24.00){\makebox(0,0)[cc]{$z_1^-$}}
\put(91.00,24.00){\makebox(0,0)[cc]{$z_1^+$}}
\put(100.33,24.00){\makebox(0,0)[cc]{$z_2^-$}}
\put(105.00,24.00){\makebox(0,0)[cc]{$z_2^+$}}
\put(115.33,24.00){\makebox(0,0)[cc]{$z_3^-$}}
\put(120.00,24.00){\makebox(0,0)[cc]{$z_3^+$}}
\put(75.33,24.00){\makebox(0,0)[cc]{$z_{0}^+$}}
\put(72.33,24.00){\makebox(0,0)[cc]{$z_{0}^-$}}
\put(59.33,24.00){\makebox(0,0)[cc]{$z_{-1}^+$}}
\put(50.67,24.00){\makebox(0,0)[cc]{$z_{-1}^-$}}
\put(40.33,24.00){\makebox(0,0)[cc]{$z_{-2}^+$}}
\put(34.67,24.00){\makebox(0,0)[cc]{$z_{-2}^-$}}
\put(26.00,24.00){\makebox(0,0)[cc]{$z_{-3}^+$}}
\put(19.50,24.00){\makebox(0,0)[cc]{$z_{-3}^-$}}
\end{picture}
\caption{ The domain $\cZ =\C\sm\cup \bar g_n$, where $g_n=(z_n^-,z_n^+)$ }
\lb{z}
\end{figure}


\begin{figure}
\tiny
\unitlength=1mm
\special{em:linewidth 0.4pt}
\linethickness{0.4pt}
\begin{picture}(120.67,34.33)
\put(20.33,20.00){\line(1,0){102.33}}
\put(71.00,7.00){\line(0,1){27.00}}
\put(70.00,18.67){\makebox(0,0)[cc]{$0$}}
\put(124.00,18.00){\makebox(0,0)[cc]{$\Re k$}}
\put(67.00,33.67){\makebox(0,0)[cc]{$\Im k$}}
\put(87.00,15.00){\linethickness{2.0pt}\line(0,1){10.}}
\put(103.00,17.00){\linethickness{2.0pt}\line(0,1){6.}}
\put(119.00,18.00){\linethickness{2.0pt}\line(0,1){4.}}
\put(71.00,15.00){\linethickness{2.0pt}\line(0,1){10.}}
\put(56.00,15.00){\linethickness{2.0pt}\line(0,1){10.}}
\put(39.00,17.00){\linethickness{2.0pt}\line(0,1){6.}}
\put(23.00,18.00){\linethickness{2.0pt}\line(0,1){4.}}
\put(85.50,18.50){\makebox(0,0)[cc]{$\pi$}}
\put(54.00,18.50){\makebox(0,0)[cc]{$-\pi$}}
\put(101.00,18.50){\makebox(0,0)[cc]{$2\pi$}}
\put(36.00,18.50){\makebox(0,0)[cc]{$-2\pi$}}
\put(117.00,18.50){\makebox(0,0)[cc]{$3\pi$}}
\put(20.00,18.50){\makebox(0,0)[cc]{$-3\pi$}}
\put(71.00,28.00){\makebox(0,0)[cc]{$ih_0$}}
\put(87.00,26.00){\makebox(0,0)[cc]{$\pi+ih_1$}}
\put(56.00,26.00){\makebox(0,0)[cc]{$-\pi+ih_1$}}
\put(103.00,24.00){\makebox(0,0)[cc]{$2\pi+ih_2$}}
\put(39.00,24.00){\makebox(0,0)[cc]{$-2\pi+ih_2$}}
\put(119.00,23.00){\makebox(0,0)[cc]{$3\pi+ih_3$}}
\put(23.00,23.00){\makebox(0,0)[cc]{$-3\pi+ih_3$}}
\end{picture}
\caption{The domain $K(h)=\C\sm\cup \G_n$, where $\G_n=(\pi n-ih_n,\pi n+ih_n)$}
\lb{k}
\end{figure}
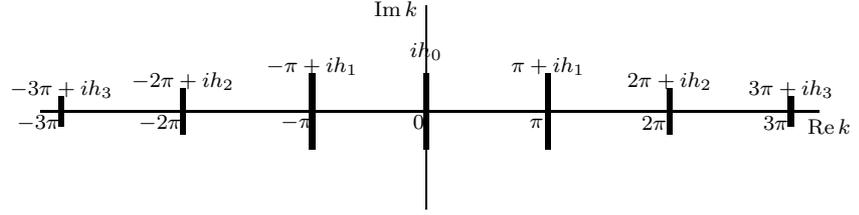
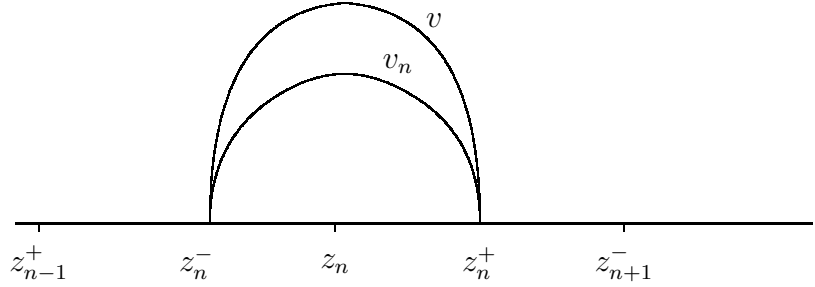
\begin{figure}
\unitlength 1mm 
\linethickness{0.4pt}
\ifx\plotpoint\undefined\newsavebox{\plotpoint}\fi 
\begin{picture}(119.75,66.5)(0,0)
\put(15,17.25){\line(1,0){104.75}}
\qbezier(40.5,17.25)(40.5,30.375)(52.25,35.75)
\qbezier(52.25,35.75)(58.25,38.375)(64.25,35.75)
\qbezier(64.25,35.75)(76.17,30.375)(76,17.25)
\qbezier(40.5,17.25)(40.5,45.125)(58.25,46.5)
\qbezier(76,17.25)(76.17,45.125)(58.25,46.5)
\put(38.75,12.00){\makebox(0,0)[cc]{$z_n^-$}}
\put(76,12){\makebox(0,0)[cc]{$z_n^+$}}
\put(65.25,38.5){\makebox(0,0)[cc]{$v_n$}}
\put(70,44.25){\makebox(0,0)[cc]{$v$}}
\put(18,17.25){\line(0,-1){1.00}}
\put(95,17.25){\line(0,-1){1.00}}
\put(57,17.25){\line(0,-1){1.00}}
\put(18,12.00){\makebox(0,0)[cc]{$z_{n-1}^+$}}
\put(95,12.00){\makebox(0,0)[cc]{$z_{n+1}^-$}}
\put(57,12.00){\makebox(0,0)[cc]{$z_{n}$}}
\end{picture}
\caption{The graph of $v(z+i0), \ z\in g_n\cup \s_n\cup \s_{n+1}$
and $|h_n|=v(z_n+i0)>0$}
\lb{grafv}
\end{figure}
%
%
%
%
%
%


Recall the following identities from \cite{KK1}
\begin{multline}
\lb{HQ}
H_0=2Q_0(h)={1\/\pi}\int\!\!\int_\C |z'(k,h)-1|^2dudv,\qqq
 H_1=4Q_1(h),\qqq
H_2=8Q_2(h),
\end{multline}
where $k=u+iv$. The functions $Q_j(h)$ are even with respect to each variable $h_n$, and then
due to \er{deav},  $Q_j(h)$ are the functions of $A=(A_n)_{n\in \Z}$.

Recall estimates for the periodic Zakharov-Shabat operator (see Corollary 2.3, \cite{K3})
\[
\lb{ae1}
{1\/2}\|q\|\leq \|h\|_2\le 3(1+\|q\|)^{1 \/2}\|q\|,
\]
\[
\lb{ae2}
{1\/2}\|\wt g\|_2\le \|q\|\le 2\|\wt g\|_2(1+\|\wt g\|_2),\qqq  \wt g=(|g_n|)_{n\in \Z},
\]
\[
\lb{ae3}
\|\e \|_2 \le 16\min \rt\{ \|q\|, \  \|h\|_2, \  \|\wt g\|_2(1+\|\wt g\|_2)   \rt\}    ,
\]
where $\e=(\e_n)_{n\in \Z}$ and $\e_n=\pi -|\s_n|\ge 0$, see \cite{K3}.

The dNLS equation admits globally defined real analytic
action-angle variables $A_n, \f_n, n\in \Z$ (see \cite{GKP}, \cite{MV1}, \cite{MV2}, \cite{VN}), where the action variables $A_n$ are given by
\[
\lb{deA}
A_n=-{1\/\pi i}\int_{\c_n}z{\D'(z)\/\sqrt{1-\D^2(z)}}dz\ge 0,
\]
see \cite{FM}, where $\c_n\ss \C$ is a counterclockwise circuit around the gap $g_n$ only and the branch of $\sqrt{1-\D^2(z)}$ is defined by $\sqrt{1-\D^2(z+i0)}>0$ for all $z\in (z_0^+,z_1^-)$.

Using the identity
$\D(z)=\cos k(z)$ and integration by parts we rewrite the action $A_n$ in terms of the quasimomentum
(see \cite{KK3})
\[
\lb{dA0}
A_n(h)=-{1\/\pi i}\int_{c_n}k(z,h)dz={2\/\pi}\int_{g_n}v(z,h)dz\ge 0, \qq v=\Im k(z), \qqq n\in \Z.
\]
Here and below the function $v(z,h)$ in the integral $\int_{g_n}v(z,h)dz$
 on the gap $g_n\ne \es$  is given by $v(z,h)=v(z+i0,h)>0, z\in g_n$.
Below we need the following estimates from [KK1]:
\[
\lb{ea}
\max\rt\{{|g_n|^2\/4}, {|g_n||h_n|\/\pi}\rt\}\le A_n={2\/\pi}\int_{g_n}v(z)dz\le {|g_n||h_n|\/\pi}
\le {2|h_n|^2\/\pi}, \qq \all \ n\in \Z.
\]
We need the important identity from \cite{KK2} (see Fig. 3)
\begin{multline}
\lb{idv1}
v(x)=v_n(x)(1+Y_n(x)),\qq Y_n(x)={1\/\pi}\int_{\R\sm g_n}{v(t)dt\/|t-x|\,v_n(t)}\ge 0,\qq \all \qq x\in g_n=(z_n^{-},z_n^{+}),\\
v_n(x)=|(x-z_n^{+})(x-z_n^{-})|^{1\/2}=|r_n^2-(z-z_n^0)^2|^{1\/2}, \qqq z_n^0={z_n^-+z_n^+\/2}, \qq r_n={|g_n|\/2}.
\end{multline}
For each $n\in \Z$ there exists the unique points $z_n\in g_n$ such that
\[
\lb{idv2}
v(z_n+i0)=|h_n|=v_n(z_n)(1+Y_n(z_n)).
\]
Define  $s=\min_{n\in \Z} |\s_n|$ and  the sequences $M_n,.., $ and the sequence $S=(S_m)_{m\in \Z}$ by
$$
M_n=\max_{z\in g_n} Y_n(z),\ \  \dot M_n=\max_{z\in g_n} |Y_n'(z)|,\ \
\ddot M_n=\max_{z\in g_n} |Y_n''(z)|,\ \
S_m={1\/2}\sum_{n\ne m}{A_n\/s^2(n-m)^2}.
$$

\begin{lemma}
\lb{esA}
Let $\|q\|\le {1\/8}$. Then for each $n\in \Z$ the following estimates hold true:
\[
\lb{pe1}
s=\min |\s_n|\ge 1,\qq \sup |g_n|\le {1\/4},
\]
\[
\lb{pe2}
\sup_{z\in g_m}{1\/\pi}\int_{g\sm g_m}
{v(t,h)dt\/(t-z_m)^2}\le S_m,
\]
\[
\lb{pe3}
M_n\le S_n,\qqq \dot  M_n\le S_n,\qqq \ddot M_n\le S_n,
\]
\[
\lb{pe4}
\|S\|_\iy\le {H_0\/2}\le {1\/128}, \qqq \|S\|_2\le {\pi^2\/6}\|A\|_2, \qqq \|S\|_1\le {\pi^2\/6}H_0\le {\pi^2\/6\cdot 64}.
\]

\end{lemma}
\no  {\bf Proof.} The estimates \er{ae3} gives $0\le \pi -|\s_n|\le 2$ and then $s\ge 1$.

The estimates \er{ae2} gives  $\sum|g_n|^2\le 4\|q\|^2\le {1\/16}$, which yields \er{pe1}.

Using $s\ge1$ and \er{A0} we get $S_m\le {1\/2}\sum_{n\ne m}A_n={H_0\/2}$. The definition of $S_m$ gives
$$
\sum_{m}S_m^2\le {1\/4} \sum_{m}\sum_{p\ne m}{1\/(p-m)^2}\sum_{n\ne m}{A_n^2\/(n-m)^2}={\pi^4\/6^2}\|A\|_2^2,
 $$
since $\sum_{n>0}{1\/n^2}={\pi^2\/6}$. Similar arguments yield the last estimate in \er{pe3} and
\er{pe2}, \er{pe4} .
\BBox

\begin{lemma}
\lb{esAh}
Let $h\in \ell^\iy$. Then for each $n\in \Z$ the following estimates hold true:
\[
\lb{esA-1}
2|h_n|\le |g_n|(1+ M_n),
\]
\[
\lb{esA-2}
0\le A_n-{|g_n|^2\/4}\le {|g_n|^2\/4}M_{n},
\]
\[
\lb{esA-4}
|h_n-\n_n|\le 4h_n M_n+h_n{|g_n|^2\/4}\rt(3(1+{|g_n|\/2})\dot M_n^2+\ddot M_n\rt),
\]
\[
\lb{pe5}
|2h_n-|g_n|(1+Y_n(z_n)|\le {|g_n|^3\/8} \dot M_n^2,
\]
\[
\lb{pe6}
|A_n-{|g_n|h_n\/4}|\le {|g_n|^4\/2^7}(\ddot M_n+6\dot M_n^2),
\]
\[
\lb{pe7}
|z_n-z_n^0|\le {|g_n|^2\/4}\dot M_n, \qqq where \qqq z_n^0={z_n^-+z_n^+\/2}.
\]
\end{lemma}
\no {\bf Proof.} Estimates \er{pe5}-\er{pe7} were proved in \cite{K4}
(see Theorem 1.3 and 1.4 in \cite{K4}).

Estimate \er{idv2} implies \er{esA-1}. Using \er{idv1} we obtain
$$
A_n={2\/\pi}\int_{g_n}v_n(x)(1+Y_n(x))dx={|g_n|^2\/4}+{2\/\pi}\int_{g_n}v_n(x)Y_n(x)dx,
$$
which gives \er{esA-2}. We will show \er{esA-4}.
Let $r={|g_n|\/2}$ and $\a_n={\n_n\/h_n}$.  We need the estimate
\[
\lb{e1}
\rt|{r^2\/h_n^2}-\a_n\rt|\le r^2C_0, \qq C_0=(3+r)\dot M_n^2+\ddot M_n,
\]
(see Lemma 2.1 in \cite{K4}). Then we obtain
\[
\lb{e2}
|\a_n-1|\le |\a_n-(r^2/h_n^2)|+|(r^2/h_n^2)-1|\le r^2C_0+2|(r/h_n)-1|.
\]
The estimate \er{pe5}
implies
\[
\rt|{r\/h_n}-1\rt|\le \rt|{r\/h_n}-{1\/1+Y_n(z_n^0)}\rt|+
{Y_n(z_n^0)\/1+Y_n(z_n^0)}\le r^3\dot M_n^2+2M_n.
\]
Then
$
|\a_n-1|\le 4M_n+r^2(3(1+r)\dot M_n^2+\ddot M_n)
$, which yields \er{esA-4}.
\BBox


\section {\bf The action $A_n(h)$ as a function of $h$}
\setcounter{equation}{0}

Define the ball $\cB^p(r)=\{\e:\|\e\|_p\le r\}\subset \ell^p, r>0$.
Let $\ell^p_{C}$ be the complexification of the space $\ell^p$.
In the complex space $\ell_{C}^p$ the corresponding ball is denoted
by $\cB_{C}^p(r)\ss \ell_{C}^p$. Define the strip
\[
\gJ_\b=\{\e \in \ell^2_{C}: \|\Im\e \|_{2} < \b\},\qqq \b>0.
\]

Recall that for fixed $k\in\C $ the function $z(k,h),\ h\in \ell^2_{\R}$ is even with respect to each variable $h_n\in \R$.
Below we need some results about the quasimomentum (see Theorem 2.3, 2.5 from  \cite{KK3}).

\begin{theorem}
\lb{Tdek}

\no i) The functional $Q_0: \ell^2\to\R_+$ has an analytic continuation
into some strip $\gJ_\b, \b>0$ and the gradient is given by:
\[
\lb{Tdek-3}
\na_n Q_0(h)=\n_n(h)=(-1)^{n-1}{\sinh h_n\/\D''(z_n)},\qqq \all \qq n\in \Z,
\qq \na_n={\pa \/ \pa h_n}.
\]

\no ii)  There exist $\ve>\b>0$ such that for any fixed real $h\in\ell^2$ the function $z(k,h+\e )$ has
the analytic extension from $(k,\e)\in\cK(h,\ve)\times \cB^2(\b)$ into the domain
$K(h,\ve)\times \cB^2_C(\b)$, where
$$
\cK(h,\ve)=\{k\in \C:\ \ \dist (k, \cup \G_n)>\ve    \}    \ss \cK(h).
$$
Moreover,  for each $n\in\Z$ the derivatives are given by (the L\"owner type equation)
\[
\lb{Tdek-1}
\na_n z(k,h)=0, \qq h_n=0,\ \ \ \ \ k\ne \pi n,  \  k\in \cK(h),
\]
\[
\lb{Tdek-2}
\na_n z(k,h)={\n_n\/z(k,h)-z_n(h)},\qq h_n\ne 0, \ \
\ \ k\ne \pi n\pm ih_n,  \  k\in\cK(h).
\]
iii) The mapping  $a:\ \ell_j^2\to \ell_j^2$ for all $j\in [0,\iy)$ given by
\[
\lb{deav}
h\to a=(a_n)_{n\in\Z }, \qqq  a_n(h)=|A_n(h)|^{1\/2}\sign h_n,
\]
is a real analytic isomorphism of $\ell_j^2=\{f=(f_n)_{n\in \Z}:
\sum_{n\in \Z} (1+n^2)^j f_n^2<\iy\}$ onto itself.

\end{theorem}

\no  {\bf Proof.} The statement (i) was proved in \cite{KK3} with $\n_n$ given by
\[
\nu_n=\ca  |k''(z_n,h)|^{-1}\sign h_n, &if \ \ \   |g_n|>0,\\
                 0,                   &if\ \ \ |g_n|=0,\ac.
\]
Using $\D(z)=\cos k$ we obtain $\D'(z)=-k'(z)\sin k$ and thus
$\D''(z_n)=-k''(z)\sin k(z_n)$, since $\D'(z_n)=0$. This and $\sin k(z_n)=i(-1)^n
\sinh h_n=$ gives $\D''(z_n)=v''(z_n)(-1)^n
\sinh h_n$, which yields \er{Tdek-1}, \er{Tdek-2}.
All other statements were proved in \cite{KK3}.
\BBox

If $h\in\ell^\iy_\R $, then for each $n\in\Z $ the following estimates from \cite{KK3} hold true:
\[
\lb{em-1}
\n_n^2\le h_n^2\le H_0,\qqq where \qqq \n_n(h)=(-1)^{n-1}{\sinh h_n\/\D''(z_n)},
\]
and recall that $z_n\in [z^-_n,z^+_n]$ is the zero of the function  $\D'(z)$.

Note that \er{Tdek-2} is the L\"owner type equation for the quasimomentum.

Introduce the contours $c_n$ around $\G_n$ and $\c_n$ around $g_n$ by
\[
\lb{decc}
c_n=\{k: {\rm dist}\ \ (k, \G_n)=\pi/4\}\ss \cK(h),\qqq \c_n=z(c_n,h)\ss\cZ.
\]
Below we need some results about the action variable $A_n$ (see Lemma 7.6 from  \cite{KK3}).

\begin{lemma}
\lb{T7.6}
Each action $A_n:\ell^2\to [0,\iy), n\in \Z$ has an analytic continuation into some strip $\gJ_\b, \b>0$ by
\[
\lb{7.11}
 A_n(h)=-{1\/\pi i}\int_{c_n}z(k, h)dk.
\]
Moreover, the derivatives ${\pa A_n\/\pa h_m}, m\in \Z$ have the following forms:
\[
\lb{7.12}
\na_m A_n(h) ={\n_m \/  \pi i}\int _{c_n}{dk\/z(k, h)-z_m },
 \qqq h\in \ell^2,\qq \na_m={\pa \/ \pa h_m},
\]
where $\n_m=\n_m(h), z_m=z_m(h)$, and if $h\in \ell^2$, then
\[
\lb{7.13}
\na_n A_n(h)=2\n_n +{2\n_n \/\pi }\int_{g\sm g_n}{v(z,h)dz \/  (z-z_n)^2},
 \qqq g=\cup g_s,
\]
\[
\lb{7.14}
\na_m A_n(h)={2\n_m \/\pi}\int_{g_n}{v'(z,h)dz\/z_m-z},
\]
\[
\lb{7.15}
\na_m A_n(h)=-{2\n_m\/\pi}\int_{g_n}{v(z,h)dz\/(z-z_m)^2},
\ \ \ m\ne n.
\]
\end{lemma}

{\bf Remark.} Recall that that ${v'(z+i0)\/z_m-z}>0$ for all $z\in g_m$ and the identities \er{7.13}, \er{7.14} and \er{pe2} give
\[
\lb{iv1+}
{1\/\pi}\int_{g_m}{v'(z)dz\/z_m-z}=1+{1\/\pi}\int_{g\sm g_m}
{v(z)dz\/(z-z_m)^2}\le 1+S_m,\qqq \all \ m\in \Z,
\]

Define the operator $F:\ell^2\to \ell^2$ by
\[
(Ff)_m=\sum_{n\in \Z} F_{m,n}f_n,\qq
\qqq F_{m,n}={\pa A_n(h)\/\pa h_m^2},\qq f=(f_n)_{n\in \Z}.
\]
Introduce the sequence $\a_n={\n_n\/2h_n}, n\in \Z$. Due to Lemma \ref{T7.6} we obtain
\[
\lb{Amn}
F_{m,n}=-{\a_m\/\pi}\int_{g_n}{v(z,h)dz \/(z-z_m)^2}, \qqq \a_m={\n_m\/h_m},\qqq
\ \ \ m\ne n,
\]
\[
\lb{Ann}
F_{n,n}=\a_n +{\a_n \/2\pi }\int_{g\sm g_n}{v(z,h)dz \/(z-z_n)^2},
 \qqq g=\cup g_s.
\]

\begin{lemma}
\lb{TdA}
Let $\|q\|\le {1\/8}$.  Then $F-I_{id}$  is the Hilbert-Schmidt operator  and  satisfies
\[
\lb{TdA-1}
|F_{m,n}|\le {A_n\/2(n-m)^2},\ \ \ \ \ m\ne n,
\]
\[
\lb{TdA-2}
|F_{n,n}-1|\le 5S_n,\qqq \qqq
\]
\[
\lb{TdA-3}
\|F-I_{id}\|_{HS}\le \pi\|q\|^2\le {\pi\/64},
\]
\[
\lb{TdA-4}
|\a_n-1|\le 5S_n,
\]
where $\|\cdot \|_{HS}$ is a Hilbert-Schmidt norm and the operator $F$ has an inverse.

\end{lemma}
\no  {\bf Proof.}   If $m\ne n$, then using \er{pe1}, \er{em-1} we obtain
$$
|F_{m,n}|\le {1\/\pi}\int_{g_n}{v(x,h)dx \/(n-m)^2}={A_n\/2(n-m)^2}.
$$
We show \er{TdA-2}. Using \er{esA-4} and Lemma \ref{esA} we obtain
$$
|\a_n-1|\le S_n\rt(4+{1\/64}\rt((3+{3\/8})S_n+1\rt)\rt)\le  {9\/2}S_n,
$$
which gives \er{TdA-4}, and
$
{\a_n \/2\pi }\int_{g\sm g_n}{v(t,h)dt \/(t-z_n)^2}\le {S_n\/2},
$
which yields \er{TdA-2}.

Combining \er{TdA-1}, \er{TdA-2} and using \er{pe4}, \er{A0}, we obtain
$$
\|F-I_{id}\|_{HS}^2\le 25\|S\|_2^2+\sum_{m}\sum_{n\ne m}{A_n^2\/4(n-m)^4}
$$
$$
\le 25\|S\|_1\|S\|_\iy+{\pi^2\/12}\|A\|_\iy^2\le
\pi^2\|q\|^4\rt({25\/48}+{1\/12}\rt)\le \pi^2\|q\|^4\le {\pi^2\/8^4},
$$
which gives \er{TdA-3}.
\BBox

\section{\bf Proof of main Theorem}
\setcounter{equation}{0}
\bigskip

Recall that the Hamiltonian has the following form (see \cite{K2})
\[
\lb{A2}
H=\sum_{n\in\Z} (2\pi n)^2A_n+2H_0^2-U,\qq \qq U(A(a))=V(h(a)), \qq
\]
where $V(h)$ is defined by
\[
\lb{AV}
V(h)={8\/3\pi}\int_{\cup g_n} v^3(z+i0,h)dz\ge0,\qqq v=\Im k(z),
\]
and $h(a)$ is an inverse mapping for $h\to a=(a_n)_{n\in\Z }, \  a_n(h)=|A_n(h)|^{1\/2}\sign h_n$ (see \er{deav})
and $A(a)$ is defined by  $A_n=a_n^2, n\in \Z$.
Recall that $a\to h$ is a real analytic isomorphism
of $\ell^2$ onto itself (see Theorem \ref{Tdek}).
In fact in order to get the properties of $U(A)$ we work with the functional $V(h)$ and below we study
how our functions depend on the Marchenko-Ostrovski parameter $h\in \ell^2$.

We define the functions $V_n(h), h\in \ell^2$, by
\[
V(h)={8\/3\pi}\int_{\cup g_n} v^3(z,h)dz=\sum_{n\in \Z}V_n, \qqq V_n(h)={8\/3\pi}\int_{g_n} v^3(z,h)dz\ge 0.
\]
Recall that $v(z+i0)> 0$ for all $z\in g_n\ne \es$ and $v$ is defined by
the equation $\cosh v(z)=(-1)^n\D(z)\ge 1, z\in g_n$. Note that for fix $k\in \cK(h)$ the function $k(z,h)$
is even with respect to each variable $h_n, n\in\Z$, and then
due to \er{A2},  $V(h(a))$ is the functions of $A=(A_n)_{n\in \Z}$.

Using $v^2(z)\le h_n^2, z\in g_n$ (see \er{idv1}) we obtain the simple estimate
\[
\lb{vn}
V_n(h)\le {8h_n^2\/3\pi}\int_{g_n} v(z+i0,h)dz={4h_n^2A_n\/3}, \qqq all \qqq n\in\Z.
\]
Consider now the function $V_n: \ell^2\to [0,\iy)$, which has good properties.

\begin{lemma}
\lb{TgH}
i) Each $V_n, n\in \Z$ has an analytic extension from $\ell^2$ into
some $\gJ_{\b}, \b>0$ by
\[
\lb{TgH-1}
V_n={4i\/\pi}\int_{c_n} (k-\pi n)^2z(k,h)dk.
\]
Moreover, if $h\in \ell^2$, then their gradients are given by
\[
\lb{TgH-2}
\na_m V_n(h)={4i \n_m(h)\/\pi}\int_{c_n} {(k-\pi n)^2\/ z(k,h)-z_m(h)}dk
={8\n_m\/\pi}\int_{g_n} {v^2(z)v'(z)\/z_m-z}dz,
\]
\[
\lb{TgH-4}
\na_m V_n(h)=-{8\n_m\/3\pi}\int_{g_n} {v^3(z)dz\/(z-z_m)^2},\qq if \qq m\ne n,
\]
\[
\lb{TgH-5}
\na_n V_n(h)={8\n_n\/\pi}\int_{g_n}{v'(z)\/z_n-z}v^2(z)dz\ge 0,
\]
where $\n_m=\n_m(h), z_m=z_m(h)$ are real analytic on $\ell^2$
and $v(z)=v(z+i0,h), z\in g_n$.

\no ii) The estimates \er{eV}, \er{TV} and the following estimates hold true:
\[
\lb{TgH-6}
|\na_n V_n(h)|\le 9|h_n|^3, \qqq and \qqq
|\na_m V_n(h)|\le {4\/3}{|h_m|h_n^2A_n\/(n-m)^2},\qq if \qq m\ne n.
\]
\end{lemma}
\no{\bf Proof.} i)
Using $v(z+i0)=-v(z-i0)$ for all $z\in g_n$ and the integration by parts,
we obtain
$$
V_n=-{4i\/3\pi}\int_{\c_n} (k-\pi n)^3dz={4i\/\pi}\int_{c_n} (k-\pi n)^2z(k)dk.
$$
The last identity and \er{Tdek-2} give
$$
\na_m V_n(h)={4i\n_m\/\pi}\int_{c_n} {(k-\pi n)^2\/z(k)-z_m}dk
={4i\n_m\/\pi}\int_{\c_n} {(k-\pi n)^2\/z-z_m}k'(z)dz,
$$
which yields \er{TgH-2}, \er{TgH-5}. Integration by parts implies \er{TgH-4}.
Note that $z_m=z_m(h)$ is real analytic on $\ell^2$, see \cite{K1}.

ii) We will show \er{TV}. Using \er{idv1} we obtain
\[
\lb{eIc}
V_n={8\/3\pi}\int_{g_n} v^3(z)dz=
{8\/3\pi}\int_{g_n} v_n^3(z)(1+Y_n(z))^3dz=I_1+I_2,
\]
$$v_n(z)=|r^2-(z-z_n^0)^2|^{1\/2},\qqq z_n^0={z_n^-+z_n^+\/2},\qqq  r={|g_n|\/2}.
$$
Now we calculate    the first integral:
\[
\lb{eII}
I_1={8\/3\pi}\int_{g_n} v_n^3(z+i0)dz={8\/3\pi}\int_{-r}^r|r^2-t^2|^{3\/2}dt=r^4.
\]
 Using \er{eIc},\er{eII} and $I_2\ge 0$ we obtain
\[
V\ge {1\/18}\sum |g_n|^4.
 \]
Using \er{esA-2} we deduce that
\[
0\le A_n^2-r^4=(A_n-r^2)(A_n+r^2)\le r^2M_n 2A_n\le 2A_n^2M_n.
\]
We estimate the second integral by
$$
0\le I_2=V_n-r^4={8\/3\pi}\int_{g_n} v_n^3(z+i0)((1+Y_n(z))^3-1)dz
$$$$
\le
M_n(1+M_n)^2{8\/\pi}\int_{g_n} v_n^3(z+i0)dz=3r^4M_n(1+M_n)^2\le
12r^4M_n,
$$
since $Y_n(z)>0$ for all $z\in g_n\ne 0$ and $M_n\le 1$ , $r^4\le A_n^2$. Thus
combine the last estimates, we get \er{TV}, since \er{pe3} and $\sum_{n>0}{1\/n^2}={\pi^2\/6}$ give
$$
\sum A_n^2M_n\le \|A\|_2^2\sup M_n\le
\|A\|_2^3 \rt(\sum_{m\ne n}{1\/(n-m)^2}\rt)^{1\/2}={\pi\/\sqrt 3} \|A\|_2^3.
$$

We show \er{TgH-6}. Using \er{TgH-4}, \er{pe1},\er{ea}, we obtain
$$
|\na_m V_n(h)|\le {8\n_mh_n^2\/3\pi}\int_{g_n} {v(z)dz\/(n-m)^2}={4\n_m\/3}{h_n^2A_n\/(n-m)^2},\qq if \qq m\ne n,
$$
and using \er{iv1+}, we get
$$
|\na_n V_n(h)|\le {8|h_n|^3\/\pi}\int_{g_n}{v'(z)dz\/z_n-z}\le 8|h_n|^3(1+S_n)\le 9|h_n|^3,
$$
which yields \er{TgH-6}.

Using \er{em-1},\er{A0} we obtain
$$
U={8\/3\pi}\int_{g} v^3(z)dz\le {8\/3\pi}\|h\|_\iy^2\int_{g} v(z)dz
={4\/3}\|h\|_\iy^2\|A\|_1\le {4\/3}\|A\|_1^2,
$$
which yields \er{eV}.
\BBox

\begin{lemma}
\lb{T32}
i) Let $q\in L^2(\T)$. Then
\[
\lb{321}
{A_n|h_n|^2\/3} \le {|g_n||h_n|^3 \/3\pi}\le V_n\le {4\/3}A_n|h_n|^2,
\]
\[
\lb{322}
{1\/3}\sum A_n|h_n|^2 \le  V\le {4\/3}\sum A_n|h_n|^2,
\]
\[
\lb{323}
\cosh h_n-1\le C_0{|g_n|^2\/8},\qqq C_0=\cosh \|h\|_\iy,
\]
\[
\lb{324}
|h_n|\le {\sqrt{C_0}\/2}|g_n|.
\]
ii) The estimate \er{eV1} holds true. Moreover, if $|h_n|=\|h\|_\iy$ for some $n\in \Z$, then
\[
\lb{325}
if \qqq C_0\ge 2, \qq \qq  \Rightarrow   \qq      |g_n|\ge 2,\qq |h_n|=\|h\|_\iy\le {\pi\/2} A_n,
\]
\[
\lb{326}
C_0\le C_1=\max \{2, \cosh {\pi\/2}\|A\|_\iy\}.
\]
\end{lemma}
\no {\bf Proof.} i) Let $g_n\ne \es$. Assume that $\a_n=z_n-z_n^-\ge |g_n|/2$, the proof of the case
$z_n^+-z_n\ge |g_n|/2$ is similar. Define the function $f_n(t)=t{|h_n|\/\a_n}, t=z-z_n^-\in (0,\a_n)$.
The function $v(z+i0), z\in g_n$ is convex and then $v(z_n^-+t+i0)\ge f_n(t), t\in (0,\a_n)$, which yields
$$
V_n={8\/3\pi}\int_{g_n} v^3(z,h)dz\ge {8\/3\pi}\int_0^{\a_n}f_n^3(t)dt={8\/3\pi}{|h_n|^3\/\a_n^3}\int_0^{\a_n}t^3dt=
{2\/3\pi}|h_n|^3\a_n\ge {1\/3\pi}|h_n|^3|g_n|.
$$
Using $v(z+i0)\le |h_n|$ for all $z\in g_n$, we get $V_n\le {8|h_n|^2\/3\pi}\int_{g_n} v(z,h)dz={4|h_n|^2\/3}A_n$.

The Taylor formula implies $\D(z_n^-)-1={1\/2}\D''(\wt z_n^-)(z_n^--z_n)^2$ for some $\wt z_n^-\in (z_n^-,z_n)$.

Using the Bernstein estimates for the bounded exponential type functions we obtain

$\sup_{z\in \R} |\D''(z)|=\sup_{z\in \R} |\D(z)|=C_0$. Then combining the Taylor formula plus the Bernstein estimates we get \er{323}, which gives \er{324}.

 ii) The estimate \er{323} yields $C_0-1\le C_0{|g_n|^2\/8},\qqq C_0=\cosh \|h\|_\iy$.
 If $C_0\ge 2$, then we deduce that ${C_0\/2}\le C_0{|g_n|^2\/8}$ and thus $|g_n|\ge 2$.
Thus, the estimate \er{ea} gives $|h_n|\le {|h_n||g_n|\/2}\le {\pi\/2}A_n$
and we obtain \er{326}.

Using \er{322},  \er{324}, \er{326}, \er{ea}, we obtain
$$
V\le  {4\/3}\sum A_n|h_n|^2 \le {2\/3}\sqrt{C_0}\sum A_n|h_n||g_n|\le {2\pi\/3}\sqrt{C_1}\|A\|_2^2.
$$
Moreover, using  \er{322}, \er{ea}, we obtain
$$
V\ge  {1\/3}\sum A_n|h_n|^2 \ge {\pi\/6}\sum A_n^2={\pi\/6}\|A\|_2^2,
$$
which yields \er{eV1}.
\BBox

The function $V(h)$ is even with respect to each variable $h_n, n\in \Z$, and then $V(h)$ is  the function of $h_n^2, n\in \Z$.

\begin{lemma}
\lb{EV}
i) Let $h\in \ell^2$ be such that $\|q\|<{1\/8}$. Then for sufficiently small
$\ve=(\ve_n)_{n\in\Z}\in\ell^2$ the following estimate holds true:
\[
|V(h+\ve)-V(h)|\le 9\|\ve\|_2.
\]
\no ii) The function $V:\ell^2\to [0,\iy)$ has the derivative $\na_n V(h)$
for each $n\in \Z$, which is continuous on $\ell^2$ and is given by
\[
\lb{31}
\na_m V(h)=\na_m V_m(h)-f_{m}(h),\qqq
f_{m}={8\n_m\/3\pi}\int_{g\sm g_m}{v^3(z)dz\/(z_m-z)^2},
\]
\[
\lb{32}
|f_{m}|\le 3|\n_m|\|A\|_\iy S_m.
\]
\end{lemma}
\no {\bf Proof.} i)  Using Lemma \ref{TgH} we get
$$
V(h+\ve)-V(h)=\sum_n (V_n(h+\ve)-V_n(h))=
\sum_{n,m\in \Z} \int_0^1\na_m V_n(h+t\ve)\ve_m dt,
$$
where all series converges absolutely.
If $\wt h=h+t\ve$, then estimates \er{TgH-6}, \er{ea} imply
$$
|\na_n V_n(\wt h)|\le 9|\wt h_n|^3,\qq and \qq
|\na_m V_n(\wt h)|\le {4|\n_m(\wt h)|\/3}{\wt h_n^2A_n(\wt h)\/(n-m)^2}
\le {8\/3\pi}{|\wt h_m \wt h_n^4|\/(n-m)^2}, \qq m\ne n,
$$
and thus using $\sum_{n>0}{1\/n^2}={\pi^2\/6}$, we get
$$
|V(h+\ve)-V(h)|\le \int_0^1\rt(\sum_{n\ne m\in \Z}{|\wt h_m|\wt h_n^4\/(n-m)^2}|\ve_n|+9\sum_{n\in \Z}|\wt h_n|^3|\ve_n|\rt)dt
$$
$$
\le \|\ve\|_2\int_0^1\rt({\pi^2\/3}\|\wt h\|_2^5+9\|\wt h\|_2^3\rt)dt\le 9\|\ve\|_2,
$$
since $\|\wt h\|\le 1$.

ii) Let $e_1=(\d_{1,n})_{n\in \Z}$ a unit vector in $\ell^p, p\ge 1$.
  Using again Lemma \ref{TgH} we obtain
$$
V(h+te_1)-V(h)=\sum_n (V_n(h+te_1)-V_n(h))=
\sum_{n} \int_0^t\na_1 V_n(h+\t e_1)d\t
$$
$$
=t\sum_{n} \na_1 V_n(h)+I_0, \qqq I_0=\sum_{n} \int_0^t\na_1 \rt(V_n(h+\t e_1)-V_n(h)\rt)d\t,
$$
where all series converges absolutely. Assume that $I_0=o(t)$ as $t\to 0$.
Then Lemma \ref{TgH} gives \er{31}.

We show that $I_0=o(t)$ as $t\to 0$. We have
$$I_0=\int_0^tG(t)d\t, \qqq G(t)=\sum_{n}  b_n(t), \qqq b_n(t)=\na_1(V_n(h+\t e_1)-V_n(h))
$$
$$
G=G_1+G_2,\qqq G_2=\sum_{-N}^N  b_n.
$$
Estimate \er{TgH-6} implies
$$
|G_1(t)|\le \sum_{|n|>N} |b_n(t)|\le \sum_{|n|>N}{C\/(1+n^2)}\le {C\/N}
$$
for some absolute  constant $C$ and $N>1$ large enough.
Each function $b_n(t,h)$ is real analytic in $h\in \ell^2$, then
$|G_2|=o(t)$ as $t\to 0$.

Similar arguments give that $\na_1 V(h)$ is continuous on $\ell^2$.

We show \er{32}. Using \er{pe2} and \er{dw-1} we obtain
$$
f_{m2}=f_{m}={8\n_m\/3\pi}\int_{g\sm g_m}{v^3(z)dz\/(z_m-z)^2}\le
{8\n_m\|h\|_\iy^2\/3\pi} \int_{g\sm g_m}{v(z)dz\/(z_m-z)^2}\le
{8\n_m\|h\|_\iy^2\/3}S_m,
$$
and \er{esA-1}, \er{pe3}, \er{pe4},  \er{esA-2} imply
$$
h_n^2\le (|g_n|^2/4)(1+S_n)^2\le A_n(1+(128)^{-1})^2, \qqq \all \ n\in \Z,
$$
which yields \er{32}.
\BBox

The function $\wt \o_n={\na_n V\/2h_n}$ is well defined and is continuous in $\ell^2$, since by Theorem \ref{Tdek}, the function $\a_n={\n_n(h)\/2h_n}$ is real analytic.

\begin{lemma}
\lb{dw}
i)  Let $\a_n={\n_n(h)\/2h_n}$. Then each component $\wt \o_m={\pa V(h)\/\pa h_m^2}, m\in \Z$
 is continuous on $\ell^2$ and satisfies for all $h\in \ell^2$:
\begin{multline}
\lb{dw-1}
\wt \o_m=\wt \o_{m1}-\wt \o_{m2},\qqq
\wt \o_{m1}={4\a_m\/\pi}\int_{g_m}v^2(z){v'(z)dz\/z_m-z}\ge 0,\qq
\wt \o_{m2}={4\a_m\/3\pi}\int_{g\sm g_m}{v^3(z)dz\/(z_m-z)^2}\ge 0.
\end{multline}
iii) If $\|q\|\le {1\/8}$, then the components $\wt \o_m, m\in \Z$ satisfy
\[
\lb{dw-2}
\wt \o_{m2}\le {9\/2}|h_n|^2, \qqq\wt \o_{m2}\le {3\|A\|_\iy\/2}S_m,
\]
\[
\lb{dw-3}
|\wt \o_{m1}-2A_m|\le 29 A_mS_m\le {A_m\/4},
\]
\[
\lb{dw-4}
|\wt \o_{m}|\le 3A_m+2\|A\|_\iy S_m,\qqq \|\wt \o\|_\iy\le 4 \|A\|_\iy,
\]
\[
\lb{dw-5}
|\wt \o_{m}-2A_m|\le 31 \|A\|_\iy S_m.
\]

\end{lemma}
\no {\bf Proof.} i) \er{TgH-4} and \er{TgH-5} imply the simple estimates \er{dw-1}.
By Lemma \ref{TgH} and Theorem \ref{Tdek}, each ${\pa V_n(h)\/\pa h_m^2}, n\in \Z$ is real analytic on $\ell^2$
and estimates \er{TgH-6} give that each component $\wt \o_m={\pa V(h)\/\pa h_m^2}, m\in \Z$
 is continuous on $\ell^2$.

ii) Using \er{32}, \er{TgH-6}, \er{em-1} we obtain \er{dw-2}.

Substituting the identity \er{iv1+} into the definition of $\wt \o_{m1}$
(see \er{dw-1}) we derive
\[
\lb{ni}
\wt \o_{m1}=4\a_m r^2+I_0-I, \qq I_0={4\a_m r^2\/\pi}\int_{g\sm g_m}
{v(t,h)dt\/(t-z_m)^2},\qq
I={4\a_m\/\pi}\int_{g_m}(r^2-v^2(t)){v'(t)dt\/(z_m-t)}.
\]
where $r={|g_m|\/2}$. Estimates \er{pe2}, \er{esA-2} give
\[
\lb{eI0}
I_0={4\a_m r^2\/\pi}\int_{g\sm g_m}
{v(t,h)dt\/(t-z_m)^2}\le 4\a_m r^2S_m\le 4A_mS_m.
\]

We rewrite  the integral $I$ in the form
\[
\lb{eI}
I={4\a_m\/\pi}\int_{g_m}(r^2-v^2(t)){v'(t)dt\/z_m-t}=I_1+I_2+I_3,\qq I_j=
{4\a_m\/\pi}\int_{g_m}f_j(t){v'(t)dt\/z_m-t},
\]
where
$$
r^2-v^2=f_1+f_2+f_3, \qq f_1=r^2-\wt v, \qq
f_2=\wt v-v_m^2, \qq f_3=v_m^2-v^2,\qq \wt v=r^2-(t-z_m)^2
$$
and recall that $v_m^2=r^2-(t-z_m^0)^2, z_m^0={z_m^-+z_m^+\/2}, t\in g_m$. We estimate all $I_j, j=1,2,3$.

Consider $I_1$. We have $f_1=(t-z_m)^2$ and then the integration by parts yields
\[
\lb{eI1}
I_1={4\a_m\/\pi}\int_{g_m}v'(t)(z_m-t)dt={4\a_m\/\pi}\int_{g_m}v(t)dt=2\a_m A_m.
\]
Consider $I_2$. Using  $f_2=(z_m-z_m^0)(2t-z_m-z_m^0)$ and \er{pe7},
\er{pe3}  we have
\[
\lb{eI2}
|f_2|\le 2|g_m||z_m-z_m^0|\le {|g_m|^3\/2}\dot M_m\le 2A_mS_m|g_m|.
\]
Consider $I_3$. Using \er{idv1} we get $f_3=v_m^2-v^2=-v_m^2Y_m(2+Y_m)$. Thus, \er{esA-2}, \er{pe3} give
\[
\lb{eI3}
|f_3|\le {|g_m|^2\/4}M_m(2+M_m)\le A_mS_m(2+S_m).
\]
Combine \er{eI}-\er{eI1} and using \er{eI2}-\er{eI3}, \er{iv1+} we have
\[
\lb{eIa}
I=2\a_m A_m+I_2+I_3,
\]
$$
|I_2+I_3|\le A_mS_mC_T{4\/\pi}\int_{g_m}{v'(t)dt\/z_m-t}
\le A_mS_mC_T(4+4S_m), \qqq
$$
where $C_T=2|g_m|+2+S_m,$
Using \er{pe1} and \er{pe4}  we obtain $C_T(4+4S_m)\le 11$, thus \er{ni}, \er{eI0}, \er{eIa} yield
\[
\lb{ewa}
\wt \o_{m1}=2\a_m (2r^2-A_m)+I_0-I_2-I_3,\qq
|I_2+I_3|\le 11A_mS_m, \qq I_0\le 4A_mS_m.
\]
Consider  $\a_m (2r^2-A_m)$, which has the form
$$
\a_m (2r^2-A_m)=A_m+I_4, \qqq I_4=(\a_m-1) A_m+2\a_m (r^2-A_m)
$$
The estimates \er{esA-2} and \er{pe2} give
\[
|r^2-A_m|\le A_mS_m.
\]
Estimate \er{TdA-4} gives
$$
|\a_m-1|\le  5S_m.
$$
Thus combine last estimates we obtain
$$
|I_4|\le 5S_mA_m+2S_mA_m=7S_mA_m,
$$
which finally together   with \er{pe4} gives
$$
\wt \o_{m1}-2A_m=2I_4+I_0-I_2-I_3, \qqq |2I_4+I_0-I_2-I_3|\le 29 S_mA_m\le A_m/4.
$$
This yields \er{dw-3}.
Combine \er{dw-2}, \er{dw-3} and  we get
\er{dw-5}, and additionally using \er{pe4} we have \er{dw-4}.
\BBox

\no {\bf Proof of Theorem \ref{T1}.}
The estimates \er{eV}, \er{TV} have been proved in Lemma \ref{TgH}.
Due to \cite{K4} the actions $A_n$ have asymptotics
$A_n=|q_n|^2(1+o(1))$ as $n\to \pm \iy$. Then if $q\in L^{4\/3}(\T)$, then
we deduce that $A\in \ell^2$ and the functional $U(A)<\iy$.
\BBox

\no {\bf Proof of Theorem \ref{T2}.}
In order to show  \er{To} we use the following identity
\[
\o_n={\pa V\/\pa A_n}=\sum_{m\in\Z}{\pa V\/\pa h_m^2}
{\pa h_m^2\/\pa A_n}=\sum_{n\in\Z}X_{n,m}\wt\o_m, \qq \wt\o_m={\pa V\/\pa h_m^2},\qq X_{n,m}={\pa h_m^2\/\pa A_n}.
\]
We rewrite the last identity in the short form:
\[
\o=X\wt\o,\qq \where \qqq \wt\o=(\wt\o_m)_{m\in \Z}, \qqq \o=(\o_m)_{m\in \Z},
\]
and $X$ is an operator in $\ell^2$ with coefficients $X_{n,m}$. Note that
$X=F^{-1}$, where the operator $F$ is defined in \er{Amn}, \er{Ann}.

Consider the operator $F$. The estimate \er{TdA-3} implies $\|F\|\ge 1-\|F-I_{id}\|\ge 1-{\pi\/64}>{15\/16}$, thus $\|F^{-1}\|\le {16\/15}$.
Let $F=I_{id}+B$. Then we obtain
$$
\o=F^{-1}\wt \o=\wt \o-F^{-1}B\wt \o=2A+f, \qqq  \where \qq f=(\wt \o-2A)-F^{-1}B\wt \o.
$$
Using \er{dw-5} we obtain
$
|\wt \o_{m}-2A_m|\le 31\|A\|_\iy S_m
$ and \er{pe4} implies
\[
\lb{z1}
\|\wt \o-2A\|_2\le 31\|A\|_\iy \|S\|_2\le    31\|A\|_\iy  {\pi^2\/6}\|A\|_2
\le 6\pi^2\|A\|_\iy \|A\|_2.
\]
Lemma \ref{TdA} yields
$$
|(B \wt \o)_m|\le |(F_{m,m}-1)\wt \o_m|+\sum_{n\ne m}|F_{m,n}\o_n|\le 5S_m|\wt \o_m|+\sum_{n\ne m}{A_n\wt \o_n\/(m-n)^2}\le 6\|\wt \o\|_\iy S_m.
$$
Then \er{pe4},\er{dw-4} give
\[
\lb{z2}
\|B \wt \o\|_2\le 6\|\wt \o\|_\iy\|S\|_2\le 4\pi^2\|A\|_\iy\|A\|_2.
\]
Combine \er{z1},\er{z2} and using $\|F^{-1}\|\le {16\/15}$ we get
$$
\|f\|_2\le \pi^2\|A\|_\iy \|A\|_2 (4+6\|F^{-1}\|)\le 11\pi^2\|A\|_\iy \|A\|_2,
$$
which yields \er{To}.
\BBox


\begin{lemma} Let $r\in [0,\pi/2]$ and $|z-\pi n|\ge r$.
Then
$$
2|\sin z|\ge e^{|\Im z|}(1-e^{-2r}).
$$

\end{lemma}

proof Sufficintely $z\in \ol\C_+$. Then
$$
2|\sin z|e^{|\Im z|}=2|e^{iz}\sin z|=|1-e^{i2z}|
$$
The max principle for the ${1\/ 1-e^{i2z}}$ for the domain
$$
\{z\in \ol\C_+: |z-\pi n|\ge r, n\in \Z\}
$$
and periodic property give that we have to check
Lemma for $|z|=r, \Im z\ge 0$. Let $w=1-e^{i2z}$.
We have
$$
|z|=r={1\/2}|\log (1-w)|\le {1\/2}\rt(|w|+{|w|^2\/2}+{|w|^3\/3}... \rt)
={1\/2}\log (1-|w|)
$$
Then $|w|\ge 1-e^{-2r}$.

\BBox

\no {\bf Acknowledgments.}
\small
The various parts of this paper were written at Mathematical Institute of the Tsukuba Univ., Japan and Ecole Polytechnique, France. The author is grateful to the Institutes for the hospitality.
I am grateful to Sergei Kuksin  for stimulating discussions and useful comments.



\begin{thebibliography}
{999}\setlength{\itemsep}{-\parskip} \footnotesize


\bibitem  [AG] {AG} Amour, L.; Guillot, J. Isospectral sets for
AKNS systems on the unit interval with generalized periodic boundary
conditions. Geom. Funct. Anal. 6 (1996), 1–27.



\bibitem  [BGGK] {BGGK} D. Battig, B. Grebert, J. Guillot, T. Kappeler,
Foliation of phase space for the cubic non-linear  Schr\"odinger
equation, Compositio Math., 85(1993), 163-199.

\bibitem  [FT] {FT} Faddeev, L. D.; Takhtajan, L. A. Hamiltonian methods in the theory of solitons. Translated from the Russian by A. G. Reyman,  Springer Series in Soviet Mathematics. Springer-Verlag, Berlin, 1987.

\bibitem  [FM] {FM} Flashka H., McLaughlin D. Canonically conjugate variables
for the Korteveg- de Vries equation and the Toda lattice with periodic
boundary conditions. Prog. of Theor. Phys. 55(1976),  438-456.



\bibitem  [GH] {GH} Gesztesy, F.; Holden, H. Soliton Equations and their algebro-geometric solutions, Vol. 1. Cambridge University Press, 2003.

\bibitem  [Go] {Go} Goluzin, G.M.  Geometric theory of functions of
a complex variable, Transl. Math. Monogr. , 26 , Amer. Math. Soc.  (1969)  (Translated from Russian).



\bibitem  [GG] {GG}  Grebert, B., Guillot, J. Gaps of one-dimensional periodic AKNS systems. Forum Math. 5 (1993), no. 5, 459--504.

\bibitem  [GKP] {GKP} Grebert, B.; Kappeler, T.; P\"oschel, P.
 Normal form theory for the NLS equation, preprint 2009.


\bibitem  [J] {J}    Jenkins A. Univalent functions and conformal mapping.
Berlin, G\"ottingen, Heidelberg, Springer, 1958.

\bibitem[KP] {KP} Kappeler, T.; P\"oschel, J. Kdv $\&$ Kam. Springer, 2003.

\bibitem[KK1] {KK1}   Kargaev, P.; Korotyaev, E. The inverse problem for the Hill
operator, a direct approach. Invent. Math. 129 (1997), no. 3, 567--593.

\bibitem[KK2] {KK2}  Kargaev, P.; Korotyaev, E. Effective masses and conformal mappings.
Comm. Math. Phys. 169 (1995), no. 3, 597--625.

\bibitem[KK3] {KK3}  Kargaev, P.; Korotyaev, E. Inverse problems generated by conformal mappings on complex plane with parallel cuts, preprint 2000,  ftp://ftp-sfb288.math.tu-berlin.de/pub/Preprints/preprint458.ps.gz

\bibitem[KK4] {KK4}  Kargaev, P.; Korotyaev, E. Identities for the Dirichlet
integral of subharmonic functions from the Cartright class. Complex
Var. Theory Appl. 50 (2005), no. 1, 35--50.

\bibitem[K1] {K1}  Korotyaev, E. Marchenko-Ostrovki mapping for periodic
Zakharov-Shabat systems. J. Differential Equations 175 (2001), no.
2, 244--274.

\bibitem[K2] {K2}   Korotyaev, E. Inverse problem and estimates for periodic
Zakharov-Shabat systems. J. Reine Angew. Math. 583 (2005), 87--115.

\bibitem[K3] {K3} Korotyaev, E. Metric properties of conformal mappings on the
complex plane with parallel cuts. Internat. Math. Res. Notices
1996, no. 10, 493--503.

\bibitem[K4] {K4} Korotyaev, E. A priori estimates for the Hill and Dirac
operators, Russ. J. Math. Phys.,15(2008), No. 3, pp. 320–-331.

\bibitem[K5] {K5} Korotyaev, E. Hamiltonian and action variables for periodic NLS,
in preparation.

\bibitem[Ku] {Ku}  Kuksin, S. B. Analysis of Hamiltonian PDEs. Oxford Lecture Series in Mathematics and its Applications, 19. Oxford University Press, Oxford, 2000.

\bibitem[Ku1] {Ku1}  Kuksin, S. private communication.

\bibitem[LS] {LS}  Levitan B., Sargsjan I. Sturm-Liouville and Dirac operators.
Translated from the Russian. Mathematics and its Applications (Soviet
Series), 59. Kluwer Academic Publishers Group, Dordericht, 1991, 350 pp.


\bibitem[L] {L} L\"owner K. Untersuchungen \"uber schlichte konforme Abbildung
des Einheitskreises, J. Math. Ann., 89(1923), 103-121.

\bibitem[LM]{LM} Lesch, M.; Malamud, M.
On the deficiency indices and self-adjointness of symmetric Hamiltonian systems,
Journal of Differential Equations
189(2003), 556--615.


\bibitem[MO]{MO} Marchenko  V.; Ostrovski I.
A characterization of the spectrum  of the Hill operator. Math. USSR  Sb.  26(1975), 493-554.

\bibitem [MV1]{MV1} McKean, H., Vaninsky, K.: Action-angle variables for the
cubic Schr\"odinger equation. Comm. Pure Appl. Math. 50 (1997),
no. 6, 489--562.

\bibitem [MV2]{MV2} McKean, H. P.; Vaninsky, K. L. Cubic Schr\"odinger: the petit canonical ensemble in   action-angle variables. Comm. Pure Appl. Math. 50 (1997), no. 7, 593--622



\bibitem[Mi1] {Mi1} Misyura T. Properties of the spectra of periodic and antiperiodic
boundary value problems generated by Dirac operators. I,II, Theor. Funktsii
Funktsional. Anal. i Prilozhen, (Russian), 30 (1978), 90-101.

\bibitem[Mi2] {Mi2} Misyura T. Properties of the spectra of periodic and antiperiodic
boundary value problems generated by Dirac operators. II, Theor. Funktsii
Funktsional. Anal. i Prilozhen, (Russian), 31(1979), 102-109.

\bibitem[V] {V} Vaninsky, K. L. Symplectic structures and volume elements in the function space for the cubic Schrödinger equation. Duke Math. J. 92 (1998), no. 2, 381--402.


\bibitem[VN] {VN} Veselov, A.; Novikov, S. Poisson brackets and complex tori. Proc. Steklov Inst. Math. 165(1985), 53–-65.




\bibitem[ZS] {ZS} Zakharov, V.E.; Shabat, A.B. A scheme for integrating nonlinear equations of mathematical physics by the method of the inverse scattering problem I. Funct. Anal. Appl. 8(1974), 226--235.





\end{thebibliography}
\end{document}